\documentclass[conference]{IEEEtran}
\IEEEoverridecommandlockouts

\usepackage{cite,color,url}
\usepackage{amsmath,amssymb,amsfonts,hyperref,multirow}
\usepackage{booktabs}
\usepackage{algorithmic,tikz}
\usepackage{subcaption,amsthm,longtable,lscape,amsmath}
\usepackage{graphicx,longtable}
\usepackage{textcomp}
\usepackage{xcolor}
\def\BibTeX{{\rm B\kern-.05em{\sc i\kern-.025em b}\kern-.08em
    T\kern-.1667em\lower.7ex\hbox{E}\kern-.125emX}}
\begin{document}

\title{On Solving the Maximum Flow \\Problem with
Conflict Constraints
}

\author{\IEEEauthorblockN{Roberto Montemanni}
\IEEEauthorblockA{\textit{Department of Sciences and Methods for Engineering} \\
\textit{University of Modena and Reggio Emilia}\\
Reggio Emilia, Italy \\
roberto.montemanni@unimore.it}
\and
\IEEEauthorblockN{Derek H. Smith}
\IEEEauthorblockA{\textit{Faculty of Computing, Engineering and Science} \\
\textit{University of South Wales}\\
Pontypridd, Wales, UK \\
derek.smith@southwales.ac.uk}}

\maketitle

\begin{abstract}
The Maximum Flow Problem with Conflict Constraints is a generalization that adds conflict constraints to a classical optimization problem on networks used to model several real-world applications. In the last few years several approaches, both heuristic and exact, have been proposed to attack the problem. In this paper we consider a  mixed integer linear program and solve it with an open-source solver. Computational results on the benchmark instances commonly used in the literature of the problem are reported. All the 160 instances benchmark instances normally used in the literature are solved to optimality for the first time, with 28 instances closed for the first time. Moreover, in the process, 6 improvements to the best-known heuristic solutions are also found.
\end{abstract}


\section{Introduction}
The study of maximum flow problems in network optimization has a long history starting in the 1950s \cite{ff56}, with research driven by real applications in the transportation
and communication fields. In this paper we consider the Maximum Flow Problem with Conflicts (MFPC), which consists of
retrieving a maximum flow from a given source to a given sink in a given network, while complying with classical capacity constraints of the connections, but also preventing positive
flow from being contemporarily carried out by the arcs of any conflicting pair. Notice that the latter restriction makes the problem considerably harder than the original version.

The MFPC was first studied in \cite{pfe13}, where the authors considered both positive disjunctive constraints (at least one of two arcs has to be selected )and negative disjunctive
constraints (at most one of two arcs can be selected), showing the NP-hardness of all the resulting variants. Later, the use of negative disjunctive
constraints to model conflicts between arc pairs in order to model real-world applications, was studied in \cite{suv20}, where some
mixed-integer linear programming formulations and exact methods 
for the problem  such as a Benders Decomposition, a Branch \& Bound and a Russian Doll Search were provided. More recently, the authors of \cite{car25} proposed a heuristic approach for the problem called Kernousel, obtained by merging two
metaheuristic approaches, one based on Carousel Greedy \cite{cer17} and the other on Kernel Search \cite{ang12}. 

Applications of the MFPC can be found in the routing of messages on telecommunication networks, where conflicts might be used to avoid the use of links of contractually-incompatible providers \cite{sch02}, in the optimization of personnel scheduling \cite{cor09}, in transportation, where conflicts constraints might be used to avoid paying alternative (cumulative) tolls while solving routing applications \cite{sch02}. Other more specific applications can be found in the image segmentation domain \cite{tar06}. 

In this paper a mixed integer linear programming model for the problem is considered and attacked by the open-source solver CP-SAT, part of the Google OR-Tools \cite{cpsat} suite for optimization. The experimental results on the benchmark instances previously adopted in the literature improve the state-of-the-art, with several improved heuristic solutions retrieved in very short time, and all the instances solved to proven optimality for the first time. This achievement is not unexpected, being in line with some recent results obtained on other combinatorial optimization problems \cite{md23}, \cite{cor}, \cite{rm25}.

The paper is organized as follows. In Section \ref{desc} the Maximum Flow Problem with Conflict Constraints is formally introduced and defined, while in Section \ref{model} the mixed integer linear programming model adopted is discussed. Section \ref{exp} is devoted to the description of the benchmark instances previously adopted in the literature and to the presentation of the experimental results obtained by the new model we propose. Finally, conclusions about the work presented are drawn in Section \ref{conc}.

\section{Problem Description}\label{desc}
The Maximum Flow Problem with additional Conflict Constraints
can be defined on a directed graph $G = (V , A)$, with $V$ denoting
the set of nodes and $A$ the set of arcs. The node set includes a source
node $s$ and a sink node $t$, while each arc $(i,j) \in A$ is associated with a capacity $u_{ij} \in \mathbb{N}^+$, indicating the maximum amount of flow that can be sent from node $i$ to node $j$ through the direct connection between them.Moreover, two arcs
$(i, j),(k,l) \in A$ are said to be \emph{in conflict} if at most one of them can have positive flow in a feasible solution. The set of arcs in conflict with the arc $(i, j) \in A$ is defined by the set $\delta(i, j)$. 
The objective of the MFPC is
to identify the maximum flow from $s$ to $t$ that satisfies the capacity and the conflict
constraints, and to route such a flow through the network.

\begin{figure*}
{
\begin{center}
{
\begin{tikzpicture}[node distance={2.5cm}, main/.style = {draw, circle}]
			\node[main,minimum size=0.75cm] (s) {s};
			\node[main,minimum size=0.75cm] (a) [above right of=s] {a};
			\node[main,minimum size=0.75cm] (b) [below right of=s] {b};
			\node[main,minimum size=0.75cm] (c) [below right of=a] {c};
			\node[main,minimum size=0.75cm] (d) [above right of=c] {d};
			\node[main,minimum size=0.75cm] (e) [below right of=c] {e};
			\node[main,minimum size=0.75cm] (t) [below right of=d] {t};
			\draw [thick, line width=1.2,->] (s) to node[midway,above left] {3} (a) ;
			\draw [thick, line width=1.2,->,color=green] (s) to node[midway,below left] {6} (b) ;
			\draw [thick, line width=1.2,->,color=blue] (a) to node[midway,below left] {1} (c) ;
			\draw [thick, line width=1.2,->,color=blue] (a) to node[midway,above] {2} (d) ;
			\draw [thick, line width=1.2,->] (b) to node[midway,left] {1} (a) ;
			\draw [thick, line width=1.2,->,color=red] (b) to node[midway,above left] {4} (c) ;
			\draw [thick, line width=1.2,->,color=blue] (b) to node[midway,below] {3} (e) ;
			\draw [thick, line width=1.2,->,color=green] (c) to node[midway,below right] {2} (d) ;
			\draw [thick, line width=1.2,->,color=red] (c) to node[midway,above right] {4} (e) ;
			\draw [thick, line width=1.2,->] (c) to node[midway,above] {3} (t) ;
			\draw [thick, line width=1.2,->] (d) to node[midway,above right] {5} (t) ;
			\draw [thick, line width=1.2,->,color=green] (e) to node[midway,below right] {7} (t) ;
\end{tikzpicture}
\hspace{1cm}
\begin{tikzpicture}[node distance={2.5cm}, main/.style = {draw, circle}]
			\node[main,minimum size=0.75cm] (s) {s};
			\node[main,minimum size=0.75cm] (a) [above right of=s] {a};
			\node[main,minimum size=0.75cm] (b) [below right of=s] {b};
			\node[main,minimum size=0.75cm] (c) [below right of=a] {c};
			\node[main,minimum size=0.75cm] (d) [above right of=c] {d};
			\node[main,minimum size=0.75cm] (e) [below right of=c] {e};
			\node[main,minimum size=0.75cm] (t) [below right of=d] {t};
			\draw [thick, line width=1.2,->] (s) to node[midway,above left] {2/3} (a) ;
			\draw [thick, line width=1.2,->,color=green] (s) to node[midway,below left] {3/6} (b) ;
			\draw [thick, line width=1.2,->,color=blue] (a) to node[midway,above] {2/2} (d) ;
			\draw [thick, line width=1.2,->,color=red] (b) to node[midway,above left] {3/4} (c) ;
			\draw [thick, line width=1.2,->] (c) to node[midway,above] {3/3} (t) ;
			\draw [thick, line width=1.2,->] (d) to node[midway,above right] {2/5} (t) ;
			\draw [thick, line width=1.2,->,color=white] (b) to node[midway,below] {3} (e) ;
\end{tikzpicture}
}
	\caption{On the left an example of a small MFPC instance is shown, where the arc capacities are placed by the arcs. Arcs in black are not affected by conflicts while any two arcs sharing the same red, blue or green color are in conflict. On the right a feasible solution carrying the maximal possible flow from $s$ to $t$ without violating any conflict is shown.}
	\label{figu}
\end{center}
}
\end{figure*}
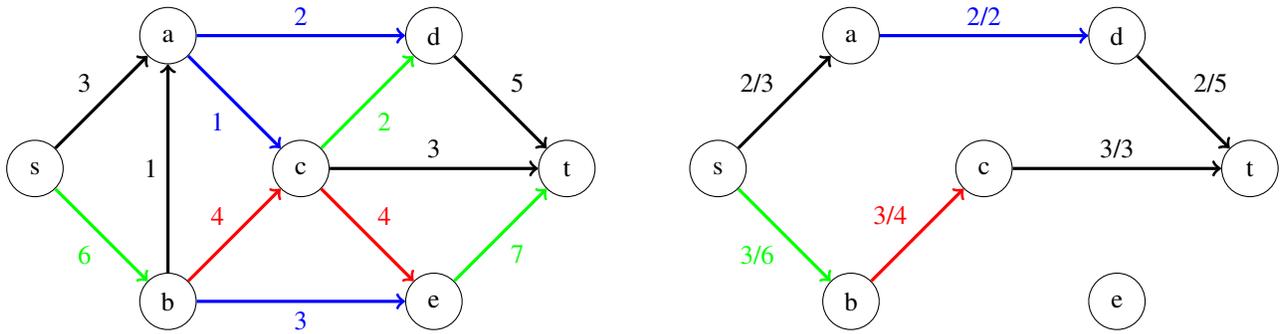

A small example of an MFPC instance (taken from \cite{car25})  and an optimal solution for the instance are provided in Figure \ref{figu}.

\section{A Mixed Integer Linear Programming Model}\label{model}
In this section a model for the MFPC \cite{suv20} is described. A variable $f_{ij} \in \mathbb{Z}^+_0$ is introduced for each $(i,j) \in A$, and represents the flow carried by the arc in the solution under investigation. In order to model conflicts, a second set of variables is introduced. For each $(i,j) \in A$, a binary variable $x_{ij}$ is considered, and takes value 1 if $f_{ij}>0$, $0$ otherwise. Finally, a variable $z \in \mathbb{Z}^+_0$ is introduced to capture the total flow from $s$ to $t$ passing through the network in the solution under investigation. The resulting model is as follows.
\begin{align} 
\ \ &  \max \ \  z& \label{1}\\ 
s.t. \ \ &	\sum_{(j,i) \in A}  \!\!\!\! f_{ji}- \!\!\!\! \sum_{(i,j) \in A} \!\!\! f_{ij} \! = \! \begin{cases}
-z  \ \text{if $i=s$}\\
+z   \ \text{if $i=r$}\\
0   \ \ \ \text{otherwise}
\end{cases} \!\!\!\!\!\!\!\!\!\!& i \in V\label{2}\\
& f_{ij} \leq u_{ij}x_{ij} & (i,j) \in A \label{3}\\
& {x_{ij}+x_{kl} \leq 1} & \!\!\!\!\!\!\!\!\!\!\!\!\!\!\!\!\!\!\!\!\!\!\!\!\!\!\!(i,j) \in A, (k,l) \in \delta_{ij} \label{4}\\
& f_{ij} \ge 0 & (i,j) \in A \label{5}\\
& x_{ij} \in \{0,1\} & (i,j) \in A \label{6}
\end{align}
The objective function (\ref{1}) maximizes the total flow from  $s$ to $t$.
Constraints (\ref{2}) deal with flow conservation and set $z$ at the same time. Notice that these constraints are the core, together with (\ref{5}), of the classic formulation for the Maximum Flow Problem. 
Inequalities (\ref{3}) regulate the activation of $x$ variables as soon as the respective $f$ variables carry some positive flow.
Constraints (\ref{4}) model conflicts and forbid conflicting edges to be active at the same time. 
Finally the domains of the variables are specified in constraints (\ref{5}) and (\ref{6}).

\section{Computational Experiments} \label{exp}
The benchmark instances adopted in the literature are introduced in Section \ref{ben}, while a comparison of the results achieved in the present work with those of the methods previously appeared in the literature, is presented in Section \ref{res}.

\subsection{Benchmark Instances}\label{ben}
The computational experiments were conducted on the collection of benchmark instances originally
proposed in \c{S}uvak et al. \cite{suv20}. By construction, each of these instances is
constructed to have at least one feasible solution with a non-zero flow,
and the characteristics of each instance depend on four parameters: 
\begin{itemize}
\item $n$, the number of nodes of the graph. Parameter $n$ takes value from the set $\{40, 50, 60, 70, 80\}$;
\item $p$. the arc density, with $p =mn(n-1)$, where $m$ denotes the number of arcs of the graph. Parameter $p$ takes value from the set $\{0.3, 0.4, 0.5, 0.6\}$;
\item $d$, the conflict density, $d= 2wm(m-1)$, with  $w$ denoting the number of conflicting
arc pairs. Parameter $d$ takes value from the set $\{0.3, 0.4, 0.5, 0.6\}$;
\item $I$, defining the arc capacity range. In setting $I=1$ the capacities are take from the interval $[10, 15]$, in setting $I=2$ the capacities are take from the interval $[15, 20]$.
\end{itemize} 

The detailed characteristics of the 160 instances obtained by combining the different possible values of the parameters can be found \cite{car25}.

\subsection{Experimental Results} \label{res}
The model discussed in Section \ref{model} has been solved with the Google OR-Tools CP-SAT solver \cite{cpsat} version 9.12. The experiments have been run on a computer equipped with an Intel Core i7 12700F CPU. The results are compared with those of the methods previously appeared in the literature. Namely, the exact methods discussed in \cite{suv20} (Benders Decomposition, \emph{BD}, Branch \& Bound, \emph{B\&B}; Russian Doll Search, \emph{RDS}) and the three metaheuristics presented in \cite{car25} (Carousel Greedy, \emph{CG}; Kernel Search, \emph{KS}; Kernousel, \emph{KO}).  The computation times of all the metaheuristic methods were obtained on a computer equipped with an Intel Core i7-3770 CPU, while the results of the exact methods were obtained on a normalized time of 2840.51 seconds on the same machine \cite{car25}. Therefore, in order to have a fair comparison, we present our computation times normalized to the same processor, according to the ratio based on the Whetstone benchmark  provided at \url{http://gene.disi.unitn.it/test/cpu_list.php}.  

The results are presented in Table \ref{t2} where for each instance (\emph{ID}) the best-known (\emph{BK}) values for lower (\emph{LB}) and upper (\emph{UB}) bounds before the present study are reported. For each exact method the percentage gap (\emph{Gap}) from the best known lower and upper bounds obtained on several runs of maximum 2840.51 seconds are presented. The gap for the lower bounds is calculated as $100 \cdot \frac{BK_{LB}-LB}{BK_{LB}}$, where $BK_{LB}$ is the best-known lower bound, and $LB$ is the value provided by the method under investigation. Analogously, the gap for the upper bound is obtained as $100 \cdot \frac{UB-BK_{UB}}{BK_{UB}}$. For each of the metaheuristic approaches we report the percentage gap for the upper bound (\emph{Gap UB}) and the computation time (\emph{Sec}). Finally, for the new approach discussed in this paper (CP-SAT) we report the optimal solution retrieved (\emph{Opt}), the time required to the retrieve the optimal solution (\emph{Sec heu}) -- useful for a comparison with the metaheuristics methods -- and the optimality-proof time (\emph{Sec opt}). Improved bounds and newly closed instances by CP-SAT are highlighted.

The results of Table \ref{t2} show that CP-SAT is able to improve the state-of-the-art results on the benchmark instances commonly adopted in the literature. In detail, all the instances are closed, with optimality proven for the first time for 28 of the 160 instances, with substantial improvements on the upper bounds. In 6 cases the value of the previously best-known heuristic solutions is improved, in some cases significantly. In terms of computation times, exact solutions are provided on average in about 200 seconds, with only two instances requiring more than 1 hour, and none above 2 hours. These computation times are already of the same order of magnitude as some of the heuristic techniques considered. When the time to retrieve the optimal solution, but without proving optimality (column \emph{Sec heu}) is considered, the time required is aligned with  (and often shorter than) that of the other heuristics, notwithstanding the higher quality of the solutions provided. When compared with the previous exact methods, the results are clearly improved.

\section{Conclusions} \label{conc}
A mixed integer linear formulation for the Maximum Flow Problem with
Conflict Constraints is considered, and solved via an open-source solver. 

The results of an experimental campaign, run on the benchmark instances adopted in the previous literature, are reported. The results indicate that the proposed approach is able to improve 6 best-known results and to find for the first time a proven optimal solution for all the 160 instances considered in the previous literature, with 28 instances closed  for the first time.

\section*{Acknowledgment}
The authors are thankful to Zeynep \c{S}uvak, Kuban Alt{\i}nel and Francesco Carrabs for sharing the benchmark instances and providing useful suggestions.

\bibliographystyle{IEEEtran}
\bibliography{mybibfile}

\onecolumn
{\scriptsize
\begin{longtable}{lrrrrrrrrrrrrrrrrrrr}
     \caption{Computational results.}\\  
     \toprule 
&&& \multicolumn{6}{c}{Previous Exact Methods}& \multicolumn{6}{c}{Previous Heuristic Methods}& \multicolumn{5}{c}{New Results}\\ 
  \cmidrule(lr){4-9}\cmidrule(lr){10-15}\cmidrule(lr){16-20}
ID  & \multicolumn{2}{c}{BK}  & \multicolumn{2}{c}{BD \cite{suv20}}& \multicolumn{2}{c}{B\&B \cite{suv20}}& \multicolumn{2}{c}{RDS \cite{suv20}} & \multicolumn{2}{c}{CG  \cite{car25}} & \multicolumn{2}{c}{KS   \cite{car25}} & \multicolumn{2}{c}{KO  \cite{car25}} & \multicolumn{5}{c}{CP-SAT}\\ 
 & & &  \multicolumn{2}{c}{Gap}&  \multicolumn{2}{c}{Gap}&  \multicolumn{2}{c}{Gap}            & Gap          & Sec          & Gap         & Sec           & Gap          & Sec        & Opt & Sec & Sec    & \multicolumn{2}{c}{Gap}  \\
    &  LB  & UB            & LB           & UB             &LB                &UB           &  LB           &    UB           & UB           &   & UB           && UB           &      &    & heu & opt    & LB &UB    \\
\cmidrule(lr){1-1}\cmidrule(lr){2-3}\cmidrule(lr){4-5}\cmidrule(lr){6-7}\cmidrule(lr){8-9}\cmidrule(lr){10-11}\cmidrule(lr){12-13}\cmidrule(lr){14-15}\cmidrule(lr){16-20}
\endfirsthead
     \caption{Computational results (continued).}\\  
     \toprule 
&&& \multicolumn{6}{c}{Previous Exact Methods}& \multicolumn{6}{c}{Previous Heuristic Methods}& \multicolumn{5}{c}{New Results}\\ 
  \cmidrule(lr){4-9}\cmidrule(lr){10-15}\cmidrule(lr){16-20}
ID  & \multicolumn{2}{c}{BK}  & \multicolumn{2}{c}{BD \cite{suv20}}& \multicolumn{2}{c}{B\&B \cite{suv20}}& \multicolumn{2}{c}{RDS \cite{suv20}} & \multicolumn{2}{c}{CG  \cite{car25}} & \multicolumn{2}{c}{KS   \cite{car25}} & \multicolumn{2}{c}{KO  \cite{car25}} & \multicolumn{5}{c}{CP-SAT}\\ 
 & & &  \multicolumn{2}{c}{Gap}&  \multicolumn{2}{c}{Gap}&  \multicolumn{2}{c}{Gap}            & Gap          & Sec          & Gap         & Sec           & Gap          & Sec        & Opt & Sec & Sec    & \multicolumn{2}{c}{Gap}  \\
    &  LB  & UB            & LB           & UB             &LB                &UB           &  LB           &    UB           & UB           &   & UB           && UB           &      &    & heu & opt    & LB &UB    \\
\cmidrule(lr){1-1}\cmidrule(lr){2-3}\cmidrule(lr){4-5}\cmidrule(lr){6-7}\cmidrule(lr){8-9}\cmidrule(lr){10-11}\cmidrule(lr){12-13}\cmidrule(lr){14-15}\cmidrule(lr){16-20}
\endhead
\multicolumn{20}{r}{\textit{(Continued on next page)}} 
\endfoot
\endlastfoot
\label{t2}
1       & 37           & 37              & 0.0                 & 0.0                  & 0.0                  & 0.0                  & 0.0                    & 0.0                      & 0.0                    & 8.0                   & 0.0                   & 4.6                  & 0.0               & 13.1              & 37     & 1.2   & 1.6    & 0.0            & 0.0            \\
2       & 49           & 49              & 0.0                 & 0.0                  & 0.0                  & 0.0                  & 0.0                    & 0.0                      & 24.5                   & 5.7                   & 0.0                   & 3.7                  & 0.0               & 9.4               & 49     & 1.3   & 1.6    & 0.0            & 0.0            \\
3       & 24           & 24              & 0.0                 & 0.0                  & 0.0                  & 0.0                  & 0.0                    & 0.0                      & 0.0                    & 1.9                   & 0.0                   & 6.3                  & 0.0               & 7.6               & 24     & 1.7   & 1.9    & 0.0            & 0.0            \\
4       & 33           & 33              & 0.0                 & 0.0                  & 0.0                  & 0.0                  & 0.0                    & 0.0                      & 0.0                    & 4.6                   & 0.0                   & 4.7                  & 0.0               & 9.6               & 33     & 1.6   & 1.9    & 0.0            & 0.0            \\
5       & 12           & 12              & 0.0                 & 0.0                  & 0.0                  & 0.0                  & 0.0                    & 0.0                      & 0.0                    & 1.5                   & 0.0                   & 12.1                 & 0.0               & 7.4               & 12     & 2.1   & 2.8    & 0.0            & 0.0            \\
\bottomrule
\end{longtable}
}
\twocolumn

\end{document}